\documentclass[A4,twoside,11pt,leqno]{article}
\usepackage{amssymb}
\usepackage{amsthm}
\usepackage[tbtags]{amsmath}
\usepackage{doc}
\usepackage{latexsym}
\usepackage{amscd}
\usepackage[frame,cmtip,arrow,matrix,line,graph,curve]{xy}
\usepackage{graphics}
\usepackage{epsfig}
\usepackage{eucal}
\usepackage{mathrsfs}
\usepackage{subfigure}
\theoremstyle{remark}

\theoremstyle{remark}
 
\renewcommand{\theequation}{\thesection.\arabic{equation}}

\font\headd=cmr8

\begin{document}
\thispagestyle{plain}
 \markboth{}{}
\small{\addtocounter{page}{0} \pagestyle{plain}
\noindent{\scriptsize KYUNGPOOK Math. J. 00(0000), 000-000}\\
\noindent{\scriptsize https://doi.org/10.5666/KMJ.0000.00.0.000}\\
\noindent {\scriptsize pISSN 1225-6951 \,\,\,\,\,\, eISSN 0454-8124}\\
\noindent {\scriptsize $\copyright$ Kyungpook Mathematical Journal}
\vspace{0.2in}\\
\noindent{\large\bf Hardy Spaces of Certain Convolution Operator }
\footnote{{}\\ \\[-0.7cm]
* Jugal K. Prajapat.
}
\vspace{0.15in}\\
\noindent{\sc Rajbala and Jugal K. Prajapat$^*$}
\newline
{\it Department of Mathematics, Central University of Rajasthan, Ajmer-305817, Rajasthan, India\\
e-mail} : {\verb|rajbalachoudary9@gmail.com|} {\it and} {\verb|jkprajapat@gmail.com|}
\vspace{0.15in}\\
{\footnotesize {\sc Abstract.} In this article, we determine conditions on the parameters of a generalized convolution operator such that it belongs to the Hardy space and to the space of bounded analytic functions. Results obtained are new and their usefulness is depicted by deducing several interesting examples.
}
\vspace{0.2in}\\
\pagestyle{myheadings}
 \markboth{\headd Rajbala and Jugal K. Prajapat$~~~~~~~~~~~~~~~~~~~~~~~~~~~~~~~~~~~~~~~~~~~~~\,$}
 {\headd Hardy Spaces of Certain Convolution Operator}\\
%
%
%
%
\noindent{\bf 1. Introduction}
\setcounter{equation}{0}
\renewcommand{\theequation}{1.\arabic{equation}}
\vspace{0.1in}\\
\indent
Let $\mathcal{H}$ denote the class of analytic functions in the open unit disk $\mathbb{D}=\{z: \; |z|<1 \},$ and $\mathcal{A}$ be the subclass of $\mathcal{H}$ consisting of normalized analytic functions of the form
$$f(z) = z + \sum_{n=2}^\infty a_n\, z^n \qquad (z \in \mathbb{D}).$$
For the analytic functions $f$ and $g$, we say that $f$ is subordinate to $g$ and write $f \prec g$, if there exists an analytic function $w$ in $\mathbb{D}$ such that $w(0)=0$ and $f(z)=g(w(z))$. In particular, if $g$ is univalent in $\mathbb{D},$ then we have the following equivalence:
$$ f(z) \prec g(z) \qquad (z \in \mathbb{D}) \quad \Longleftrightarrow \quad [f(0) = g(0) \quad \mbox{and} \quad f(\mathbb{D}) \subset g(\mathbb{D})]. $$
Further, for functions $f_j \in \mathcal{A}$, given by \,$f_j(z) = z + \sum_{n=2}^{\infty} a_{n, j} \; z^n \;\; (j=1,2),$ the Hadamard product (or Convolution) of $f_1$ and $f_2$ is defined by
$$(f_1*f_2)(z) :=  z+ \sum_{n=2}^{\infty} a_{n,1} \;a_{n,2} \, z^{n} \qquad (z \in \mathbb{D}). $$

Let $\mathcal{R}(\beta)$ denote the class of functions $f \in \mathcal{A}$ such that $\Re(f'(z))>\beta \;\; (z \in \mathbb{D}, \,\beta<1),$ and let $\mathcal{P}(\beta)$ denote the class of functions $f \in \mathcal{H}$ such that $f(0)=1$ and \,$\Re(f(z))>\beta \;\;(z \in \mathbb{D}, \,\beta<1).$ For $\beta=0$, we denote $\mathcal{R}(\beta)$ and $\mathcal{P}(\beta)$ simply by $\mathcal{R}$ and $\mathcal{P}$, respectively. 

Let $\mathcal{H}^\infty$ denote the space of all bounded analytic functions in $\mathbb{D}$. This is Banach algebra with respect to the norm $\Vert f\Vert_{\infty}=\mbox{sup}_{z \in \mathbb{D}} \,|f(z)|$. For $f \in \mathcal{H}$, we define
\begin{equation}\label{hardy1}
M_p(r, f)=\left\{\begin{array}{ll}
\displaystyle\left(\frac{1}{2\pi}\int_0^{2\pi} \left|f(re^{i\theta})\right|^p \,d\theta \right)^{1/p}  & \; (0<p<\infty),\\
\displaystyle\max_{|z|\leq r} |f(z)| & \; (p=\infty).
\end{array}\right.
\end{equation}
The function $f \in \mathcal{H}$ belongs to the Hardy space $\mathcal{H}^p\,(0<p\leq \infty),$ if  $M_p(r, f)$ is bounded for all $r \in [0, 1)$. Clearly, we have 
$$\mathcal{H}^\infty \subset \mathcal{H}^q \subset \mathcal{H}^p \quad \mbox{for} \quad 0<p<q<\infty $$
(see \cite[p. 2]{PLD}). For $1\leq p \leq \infty, \; \mathcal{H}^p$ is a Banach space with the norm defined by 
\begin{equation}\label{norm}
\Vert f \Vert _p = \lim_{r \rightarrow 1^-} {M_p(r, f)}\quad \quad (1\leq p \leq \infty)
\end{equation}
(see \cite[p. 23]{PLD}). Following are two widely known results (see \cite{YK}) for the Hardy space $\mathcal{H}^p$: 
\begin{eqnarray}\label{hardy2}
\Re(f'(z))>0 &\Longrightarrow & f' \in \mathcal{H}^p \;\; \mbox{for \,all} \;\;p<1 \\
&\Longrightarrow & f \in \mathcal{H}^{q/(1-q)} \;\; \mbox{for \,all} \;\;0<q<1. \notag
\end{eqnarray}

In \cite{samy}, Ponnusamy studied the Hardy space of hypergeometric functions. Further, Baricz \cite{baricz} obtained the conditions for the generalized Bessel functions such that it belongs to Hardy space and Yagmur and Orhan \cite{yagmur} studied the same problem for the generalized Struve functions. 
 
In \cite{IROG} (see also \cite{RWIS}), Ibrahim  studied the following generalized fractional integral operator in the complex plane $\mathbb{C}$: 
\begin{equation}\label{genop}
I^{\alpha, \mu}_zf(z)=\frac{(\mu+1)^{1-\alpha}}{\Gamma(\alpha)}\int^z_0(z^{\mu+1}-\xi^{\mu+1})^{\alpha-1} \xi^{\mu}f(\xi)d\xi\quad\quad (\alpha, \mu  \in\mathbb{R},\alpha>0,  \mu \neq -1 ),
 \end{equation}
where the function $f(z)$ is analytic in a simply connected region of $\mathbb{C}$ containing the origin, and the multiplicity of $(z^{\mu+1}-\xi^{\mu+1})^{\alpha-1}$ is removed by requiring $\log(z^{\mu+1}-\xi^{\mu+1})$ to be real when $(z^{\mu+1}-\xi^{\mu+1})>0.$ We observe that, if we take $\mu=0$ in \eqref{genop}, we arrive at the {\it Srivastva-Owa}  fractional integral operator \cite{HMSS} (see also \cite{JKP}) and when $\mu\rightarrow -1^+$ in \eqref{genop}, we obtain the {\it Hadamard} fractional integral operator \cite{HJE}. 

For $\alpha>0$ and $\mu> -1$, we define a fractional integral operator $\Omega^{\alpha, \mu}_z:\mathcal{A} \rightarrow \mathcal{A}$ by
\begin{eqnarray}
\Omega^{\alpha,\mu}_z f(z) &=&  \dfrac{(\mu+1)^{\alpha}\Gamma\left(\frac{1}{\mu+1}+\alpha+1\right)}{\Gamma\left(\frac{1}{\mu+1}+1\right)}z^{-\alpha(\mu+1)} I^{\alpha, \mu}_zf(z) \notag \\
&=& z+\sum_{n= 2}^{\infty}\dfrac{\Gamma\left(\frac{1}{\mu+1}+\alpha+1\right) \Gamma\left(\frac{n}{\mu+1}+1\right)}{\Gamma\left(\frac{n}{\mu+1}+\alpha+1\right) \Gamma\left(\frac{1}{\mu+1}+1\right)}\, a_n z^n.  
\end{eqnarray} 
Note that\; $\Omega^{0, 0}_zf(z)= f(z), \; \;\Omega^{1, 0}_zf(z)=\mathcal{L}f(z)$  and $\Omega^{\alpha, 0}_zf(z)=L(2, \alpha+2)f(z),$ where $\mathcal{L}$ and $L(a,c)$ denotes the Libera integral operator \cite{Libera} and  Carlson-Shaffer operator \cite{BCCD}, respectively. Also, we observe that the operator $\Omega _z ^{\alpha, \mu}$ satisfies the following recurrence relation 
 \begin{equation}\label{eqn4}
 z \,(\Omega^{\alpha,\mu}_zf(z))'=(1+\alpha(\mu+1))\Omega^{\alpha-1,\mu}_z f(z)-\alpha(\mu+1)\Omega^{\alpha,\mu}_z f(z) \quad (\alpha>1, \,\mu> -1).
\end{equation} 
Corresponding to fractional integral operator $\Omega_z^{\alpha, \mu}f(z),$ a fractional differential operator $\Phi^{\alpha, \mu}f(z)$ was studied in \cite{IROG, RWIS} to obtain its boundedness in Bergman space and certain geometric properties were also discussed.
\vspace{0.2in}\\
%
%
%
%
%
\setcounter{equation}{0}
\renewcommand{\theequation}{2.\arabic{equation}}
\noindent{\bf 2. Main Results}
\vspace{0.1in}\\
\setcounter{equation}{0}

In order to derive our main results, we recall here the following lemmas:

\medskip
\noindent
{\bf Lemma 2.1.} \;{\rm(Suffridge \cite{TSS})} \label{lema1}
Let $F$ and $G$ be analytic functions in $\mathbb{D}$ and $F(0)= G(0).$ If $H(z)= zG'(z)$ is starlike function in $\mathbb{D}$ and $zF'(z)\prec zG'(z)$, then
$$F(z)\prec G(z)= G(0)+\int^z_0\frac{H(t)}{t}dt.$$
 
\noindent
{\bf Lemma 2.2} \;{\rm(Stankiewicz and Stankiewicz \cite{stank})}\label{lema2}
For $\alpha <1$, $\beta <1$, we have $\mathcal{P}(\alpha) * \mathcal{P}(\beta) \subset \mathcal{P}(\delta),$ where $\delta=1-2(1-\alpha)(1-\beta)$. The value of $\delta$ is best possible.

\medskip
\noindent
Our first main result is given by Theorem 2.3 below.
\medskip

\noindent
{\bf Theorem 2.3}\label{thm2}
{\it Let $\alpha >1$ and $\mu > -1$. If $f(z)\in \mathcal{R}$, then $\,\Omega_z^{\alpha, \mu}f(z) \in \mathcal{R} \cap\mathcal{H}^{\infty}$.} 
\vspace{0.05in}\\
\noindent
{\it Proof.} From the definition of operator $\Omega_z^{\alpha,\mu}f(z)$, we have 
\begin{equation}\label{13}
 \Re \left\{(\Omega_z^{\alpha,\mu}f(z))'\right\} = \dfrac{(\mu+1) \Gamma\left(\frac{1}{\mu+1}+\alpha+1\right)}{\Gamma(\alpha) \Gamma\left(\frac{1}{\mu+1}+1\right)}\int_0^1 u^{\mu}(1-u^{\mu+1})^{\alpha-1} \,\Re\{f'(zu)\} \,du.
 \end{equation}
By hypothesis $f(z) \in \mathcal{R}$, hence it follows that $\Omega_z^{\alpha, \mu}f(z) \in \mathcal{R}.$ By first implication of \eqref{hardy2}, we have $\left(\Omega_z^{\alpha, \mu}f(z)\right)' \in \mathcal{H}^q$ for all $q<1$. Further, by second implication of \eqref{hardy2}, we have $\Omega_z^{\alpha, \mu}f \in \mathcal{H}^{q/(1-q)}$ for all $0<q<1$, or equivalently, $\Omega_z^{\alpha, \mu}f \in \mathcal{H}^p$ for all $0<p<\infty$.

Since $f \in \mathcal{R}$, then using the well known bound for {\it Caratheodory functions} in $\mathbb{D}$, we have $|a_n|\leq 2/n, \; n \geq 2$. Hence
\begin{eqnarray}\label{11}
\left|\Omega_z^{\alpha, \mu}f(z) \right| & \leq & |z|+\frac{\Gamma(\frac{1}{\mu+1}+\alpha+1)}{\Gamma(\frac{1}{\mu+1}+1)}\sum_{n=2}^{\infty} \frac{\Gamma(\frac{n}{\mu+1}+1)}{\Gamma(\frac{n}{\mu+1}+\alpha+1)} |a_n| |z|^n\\
&<& 1+2 \sum_{n=2}^{\infty}\theta(n), \notag
\end{eqnarray}
where $\theta(n)$ is given by
 \begin{equation}{\label{P}}
\theta(n)= \dfrac{\Gamma(\frac{1}{\mu+1}+\alpha+1)\,\Gamma(\frac{n}{\mu+1}+1)}{n \,\Gamma(\frac{n}{\mu+1}+\alpha+1)\,\Gamma(\frac{1}{\mu+1}+1)}.
\end{equation}
Evidently, $\theta(1)=1$ and $\theta(n)>0$ for all $n \in \mathbb{N}.$ It is well known that the digamma function $\Psi(x)=\Gamma'(x)/\Gamma(x)$ is increasing for all $x>0$ (see \cite[II, sec. 3/ eqn. 4 on p. 723]{APYB}). Therefore $\Psi(x+\varepsilon)> \Psi(x)$\, for $\varepsilon>0.$ The auxiliary function $\tilde{\Gamma}(x):=\Gamma(x+\varepsilon)/\Gamma(x)\,$ has a positive derivative as
\begin{eqnarray}
\tilde{\Gamma}'(x) &= &\dfrac{\Gamma'(x+\varepsilon)\Gamma(x)-\Gamma(x+\varepsilon)\Gamma'(x)}{[\Gamma(x)]^2} \notag \\
&=& \dfrac{\Gamma(x+\varepsilon)}{\Gamma(x)}\left(\dfrac{\Gamma'(x+\varepsilon)}{\Gamma(x+\varepsilon)}-\dfrac{\Gamma'(x)}{\Gamma(x)}\right)>0 \quad (x>0, \;\varepsilon>0).
\end{eqnarray}
Thus $\tilde{\Gamma}(x)$ is an increasing function, so that
\begin{equation}\label{epsilon}
\dfrac{\Gamma(x+\varepsilon)}{\Gamma(x)}\geq \dfrac{\Gamma(y+\varepsilon)}{\Gamma(y)} \quad \mbox{whenever} \quad x\geq y>0.
\end{equation}
Hence
$$\dfrac{\theta(n)}{\theta(n+1)}= \dfrac{(n+1)\,\Gamma(1+\frac{n}{\mu+1})\,\Gamma(1+\alpha+\frac{n+1}{\mu+1})}{n \,\Gamma(1+\frac{n+1}{\mu+1})\,\Gamma(1+\alpha+\frac{n}{\mu+1})} \geq 1, $$
which shows that $0< \theta(n+1)<\theta(n)\leq \theta(2)$\, for each $\,n \geq 2.$ Hence $\theta(n)$ is a non-increasing function of $n \geq 2.$ 

Now, we shall show that $\lim _{n\rightarrow\infty }|\theta(n)|^{1/n}=1.$ Using the asymptotic formula for the Gamma function \cite[sec. 1.18, (4)]{AEWMF}, we have 
\begin{eqnarray*}
\lim_{n\rightarrow \infty}|\theta(n)|^{1/n}&=& \lim_{n \rightarrow \infty} \left[\dfrac{\Gamma(1+\alpha+\frac{1}{\mu+1})}{\Gamma(1+\frac{1}{\mu+1})} \right]^{1/n}.\lim_{n \rightarrow \infty}\left[\dfrac{\Gamma(1+\frac{n}{\mu+1})}{n \,\Gamma(1+\alpha+\frac{n}{\mu+1})} \right]^{1/n} \notag \\
& \sim & \lim_{n \rightarrow \infty}\left[\left(\frac{n}{\mu+1}\right)^{-\alpha}\right]^{1/n}=\lim _{n \rightarrow \infty}(n^{1/n})^{-\alpha}\left[\frac{1}{(\mu+1)^{\alpha}}\right]^{1/n}=1.
\end{eqnarray*}
This implies that the series in \eqref{11} converges for $|z|<1$. Further, applying Raabe's test for convergence, we deduce that the power series for $\Omega_z^{\alpha, \mu}f$ converges absolutely for $|z|=1$ for all $\alpha>0$ and $\mu>-1$.

Also, using the known result \cite[Theorem 3.11]{PLD}, $\left(\Omega_z^{\alpha, \mu}f \right)' \in \mathcal{H}^q$ implies continuity of $\Omega_z^{\alpha, \mu}f$ on the compact set $\overline{\mathbb{D}}$. Since the continuous function $\Omega_z^{\alpha,\mu}f$ on the compact set $\,\overline{\mathbb{D}}$ is bounded, thus $\Omega_z^{\alpha,\mu}f$ is a bounded analytic function in $\mathbb{D}$. Therefore $\Omega_z^{\alpha,\mu}f  \in \mathcal{H}^{\infty}$ and this completes the proof.
\hfill$\Box$
\medskip

\noindent 
{\bf Remark 2.4}. If we consider $\Psi(z)=-z-2\log (1-z)=z+2\sum_{n=2}^\infty \dfrac{z^n}{n}$, then from the geometrical descriptions of image domains of functions in the class $\mathcal{R},$ we can easily verify that $\Psi \in \mathcal{R}$ and not belongs to $\mathcal{H}^\infty,$ but as per the Theorem 2.3, function $\,\Omega_z^{\alpha, \mu}\Psi(z) \in \mathcal{H}^{\infty} \cap \mathcal{R}$, where
\[\Omega^{\alpha,\mu}_z \Psi(z) = z+2\sum_{n= 2}^{\infty}\dfrac{\Gamma\left(\frac{1}{\mu+1}+\alpha+1\right) \Gamma\left(\frac{n}{\mu+1}+1\right)}{n \; \Gamma\left(\frac{n}{\mu+1}+\alpha+1\right) \Gamma\left(\frac{1}{\mu+1}+1\right)} \, z^n   \quad (\alpha >1, \,\mu > -1, z \in \mathbb{D}).\]

\medskip
\noindent
{\bf Theorem 2.5}  Let $f \in \mathcal{A}, \alpha >1$ and $\mu>-1.$ If 
\begin{equation}\label{corollary}
\left|\dfrac{\Omega^{\alpha-1,\mu}_zf(z)- \Omega^{\alpha, \mu}_zf(z)}{z}\right|< \dfrac{1}{1+\alpha(\mu+1)}\qquad (z\in \mathbb{D}),
\end{equation}
then \;$\dfrac{\Omega^{\alpha, \mu} _z f(z)}{z} \in \mathcal{P}.$ 
 
{\it Proof.} The inequality \eqref{corollary} is equivalent to 
$$\dfrac{\Omega^{\alpha-1, \mu}_zf(z)- \Omega ^{\alpha,\mu}_z f(z)}{z}\prec \dfrac{z}{1+\alpha(\mu+1)} \qquad (z \in \mathbb{D}),$$
which in view of  \eqref{eqn4}, can be written as $z\left(\dfrac{\Omega ^{\alpha, \mu}_z f(z)}{z}\right)' \prec z(1+z)'.$ Now applying Lemma 2.1,  for $F(z)= \dfrac{\Omega^{\alpha, \mu}_z f(z)}{z}$\, and\, $G(z)=1+z,$ the desired result follows.

\medskip
\noindent
{\bf Theorem 2.6.} Let $f \in \mathcal{A}, \alpha > 1$ and $\mu >-1.$ If $g\in \mathcal{R}(1/2)$ and $\Omega_z ^{\alpha, \mu} f$ satisfy the inequality (\ref{corollary}), then $\Omega_z ^{\alpha, \mu} f (z)* g(z)\in \mathcal{R}.$
 
{\it Proof.} Let $u(z)= \Omega_z ^{\alpha, \mu}f(z)*g(z),$ then $u'(z)=\dfrac{\Omega_z ^{\alpha, \mu} f(z)}{z}*g'(z).$ In view of the hypotheses and Theorem 2.5, we have $\dfrac{\Omega_z ^{\alpha, \mu} f(z)}{z}\in \mathcal{P}\; \; \mbox{and}\; \; g'(z) \in \mathcal{P}(1/2).$ Now using Lemma 2.2, we obtain that $u'(z) \in \mathcal{P},$ which is equivalent to $u(z)\in \mathcal{R},$ hence the desired result follows.

\medskip
{\bf Acknowledgement}. The authors wish to thank the referee for her/his valuable suggestions, which have improved the paper. The work was supported by the University Grant Commission, India (Project Number: MRP-MAJOR-MATH-2013-19114).

--
\vspace{0.1in}\\

\footnotesize{
}

\end{document}